\documentclass[12pt]{article}

\usepackage{amsmath, epsfig, cite}
\usepackage{amssymb}
\usepackage{amsfonts}
\usepackage{latexsym}
\usepackage{float}
\usepackage{color}
\usepackage{mathrsfs}

\newtheorem{thm}{Theorem}[section]

\newtheorem{lem}[thm]{Lemma}

\newcommand{\pf}{\noindent{\it Proof} }

\numberwithin{equation}{section}

\newcommand{\qed}{{\hfill$\square$}\medskip}

\begin{document}
	\begin{center}
		{\large\bf  Proof of Two Supercongruences of Guillera and Zudilin      }
	\end{center}
	\vskip 2mm \centerline{Wei-Wei Qi}
	
	\begin{center}
		{\footnotesize MOE-LCSM, School of Mathematics and Statistics, Hunan Normal University, Hunan 410081, P.R. China\\[5pt]
			{\tt wwqi2022@foxmail.com} \\[10pt]
		}
	\end{center}
	
	\vskip 0.7cm \noindent{\bf Abstract.} In $2012$,  Guillera and Zudilin established the following two supercongruences involving truncated Ramanujan-type series: for any odd prime $p>2$, 
	\begin{align*}
		\sum_{n=0}^{p-1}\frac{(\frac{1}{2})_n(\frac{1}{3})_n(\frac{1}{4})_n(\frac{3}{4})_n}{(1)_n^5}(-1)^n\left(172n^2+75n+9\right)\left(\frac{27}{16}\right)^n\equiv 9p^2 \pmod{p^5},
	\end{align*}	
and
\begin{align*}
	\sum_{n=0}^{p-1}\frac{(\frac{1}{2})_n(\frac{1}{3})_n(\frac{2}{3})_n}{(1)_n^3}\left(11n+3\right)\left(\frac{27}{16}\right)^n\equiv 3p \pmod{p^3},
\end{align*} 
	where $(a)_n=\prod_{k=0}^{n-1}(a+k)$ denotes the Pochhammer symbol (rising factorial). In this paper, we mainly apply the Wilf-Zeilberger (WZ) method and symbolic summation techniques to prove these two supercongruences.

	\vskip 2mm
	\noindent{\it MR Subject Classifications}: 11A07, 05A10, 11B65, 33C20.	
	
	\section{Introduction} 
	
	In $1914$, Ramanujan \cite{wt-0-1} presented a number of infinite series for $1/\pi$  in the form
	\begin{align*}
		\sum_{k=0}^{\infty}(ak+b)\frac{t_k}{m^k}=\frac{c}{\pi},
	\end{align*}		
	which $a,b,c,m$ are constants and $t_k$ denotes a product of the binomial coefficients ${2k\choose k}$, ${3k\choose k}$, ${4k \choose 2k}$, ${6k\choose 3k}$. In $2020$, Guillera \cite{wt-0-2} explain a general method for proving Remanaujan-type series for $1/\pi$.
	In $1997$, Van Hamme \cite{wt-0-3} investigated the partial sums of Ramanujan-type series for $1/\pi$, and put forward thirteen conjectured supercongruences labeled as $(A.2)-(M.2)$, which take the following form:
	 \begin{align*}
	 	\sum_{k=0}^{\frac{p-1}{2}}(ak+b)\frac{t_k}{m^k}=bp\left(\frac{\varepsilon_d d}{p}\right) \pmod{p^r},
	 \end{align*}
	where $d$ is a positive square-free number, $r$ is a positive integer, $p$ is an odd prime satisfying $p\nmid dm$, $\varepsilon_d \in \{\pm1\}$ with $\varepsilon_d=-1$ if $d>1$, $\left(\frac{\cdot}{p}\right)$ denotes the Legendre symbol. This type of congruences is particular intriguing and has been investigated by many scholars.  In $1990$, H. S. Wilf and D. Zeilberger \cite{wt-0-4} introduced the WZ method, which can  discover and prove hypergeometric identities. In $1995$, Gessel \cite{wt-0-5} developed a systematic approach to construct WZ pairs. Recently, Feng and Hou (\cite{wt-3}, \cite{wt-4})  investigated several new congruences associated with truncated Ramanujan-type series via the WZ method. For more applications of WZ method, readers may refer to (\cite{wt-0-6}, \cite{wt-0-7}, \cite{wt-0-8}, \cite{wt-0-9}, \cite{wt-0-10} ) and other relevant literature.

In $2012$,	Guillera and Zudilin \cite[(10) and (12)]{wt-1} showed that
\begin{align*}
	\sum_{n=0}^{\frac{p-1}{2}}\frac{\left(\frac{1}{2}\right)_n^3}{n!^3}(3n+1)2^{2n}\equiv p \pmod{p^3} \quad for \quad p>2,
\end{align*}	
\begin{align*}
	\sum_{n=0}^{p-1}\frac{\left(\frac{1}{2}\right)_n^3}{(1)_n^3}(3n+1)(-1)^n2^{3n}\equiv (-1)^{\frac{p-1}{2}}p \pmod{p^3} \quad for \quad p>2.
\end{align*}
In their work, the authors proposed the following two supercongruence conjectures \cite[(7) and (9)]{wt-1}: for any prime $p>2$
\begin{align}
	\sum_{n=0}^{p-1}\frac{(\frac{1}{2})_n(\frac{1}{3})_n(\frac{1}{4})_n(\frac{3}{4})_n}{(1)_n^5}(-1)^n\left(172n^2+75n+9\right)\left(\frac{27}{16}\right)^n\equiv 9p^2 \pmod{p^5}, \label{con-2}
\end{align}		
and	
\begin{align}
	\sum_{n=0}^{p-1}\frac{(\frac{1}{2})_n(\frac{1}{3})_n(\frac{2}{3})_n}{(1)_n^3}\left(11n+3\right)\left(\frac{27}{16}\right)^n\equiv 3p \pmod{p^3}. \label{con-1}
\end{align}
And they mentioned that although they had also obtained the WZ pairs for  \eqref{con-2} and \eqref{con-1}, these pairs did not seem to be sufficient to prove that the corresponding congruences hold modulo the expected powers of $p$. Recently, Feng and Hou \cite[Theorem $1.4$ and Theorem $1.1$]{wt-4} constructed two distinct the WZ pairs, by means of which they proved the supercongruence \eqref{con-2} and \eqref{con-1} modulo $p^3$ and $p^2$, respectively.
		
Motivated by the above work, the main goal of this paper is to prove the above  aforementioned two supercongruence conjectures.  
The Fermat quotient of an integer $a$ with respect to an odd prime $p$ is gives by
\begin{align*}
	q_p(a)=\frac{a^{p-1}-1}{p}.
\end{align*}
Obviously, for nonnegative integer $n$,	
	 \begin{align*}
		\frac{(\frac{1}{2})_n(\frac{1}{3})_n(\frac{1}{4})_n(\frac{3}{4})_n}{(1)_n^5}(-1)^n\left(172n^2+75n+9\right)\left(\frac{27}{16}\right)^n=\frac{(-1)^n(172n^2+75n+9)}{2^{12n}}{2n \choose n}^3{3n \choose n}{4n\choose 2n},
	\end{align*}	
	and	
	\begin{align*}
		\frac{(\frac{1}{2})_n(\frac{1}{3})_n(\frac{2}{3})_n}{(1)_n^3}\left(11n+3\right)\left(\frac{27}{16}\right)^n=\frac{11n+3}{64^n}{2n \choose n}^2{3n \choose n}.
	\end{align*}	
\begin{thm}
	For any odd prime $p>2$, we have
	\begin{align*}
		\sum_{k=0}^{\frac{p-1}{2}}\frac{(-1)^k(172k^2+75k+9)}{2^{12k}}{2k \choose k}^3{3k \choose k}{4k\choose 2k}
		\equiv 9p^2+6p^3q_p(2)-9p^4q_p(2)^2 \pmod{p^5}. 
	\end{align*}
\end{thm}	
\begin{thm}
	\eqref{con-2} is true.
	\end{thm}
	\begin{thm}
	\eqref{con-1} is true.
\end{thm}		
\begin{thm}
		For any odd prime $p>2$, we have
	\begin{align*}
		\begin{aligned}
			\sum_{k=0}^{\frac{p-1}{2}}\frac{11k+3}{64^{k}}{2k \choose k}^2{3k \choose k}
			\equiv 3p+3p^2q_p(2) \pmod{p^3}. 
		\end{aligned}	
	\end{align*}	
\end{thm}

	The structure of this paper is organized as follows. In Section $2$, we present some necessary Lemmas. We are going to prove Theorems $1.1$,  $1.2$,  $1.3$, and $1.4$ in Section $3$,  $4$, $5$ and  $6$, respectively.	
	
\section{Some Lemmas}
	Let $H_n(r)$ denote the $n$-th generalized harmonic number of order $r$, which is defined as
	\begin{align*}
		H_n(r)=\sum_{k=1}^{n}\frac{1}{k^r},
	\end{align*}
	with the convention that $H_n=H_n(1)$. Additionally, we adopt the following notation:	
	\begin{align*}
		H_n(-r)=\sum_{k=1}^{n}\frac{(-1)^k}{k^r} \quad and \quad H(1,1;n)=\sum_{1\leq i<j\leq n}\frac{1}{ij}.
	\end{align*}	
	
Recall the classical Wolstenholme congruence \cite{wt-5} says that	
\begin{align}
	H_{p-1}\equiv 0 \pmod{p^2}, \label{wl0-1}
\end{align}
\begin{align}
	H_{p-1}(2)\equiv 0 \pmod{p}. \label{wl0-2}
\end{align}

\begin{lem} For any prime $p>2$, we have			
		\begin{align}
			&\sum_{k=1}^{\frac{p-1}{2}}\frac{H_{2k}}{k}\equiv	q_p(2)^2 \pmod{p}, \label{wl0-3} \\ 
				&\sum_{k=1}^{\frac{p-1}{2}}\frac{1}{k}\equiv -2q_p(2)+q_p(2)^2p   \pmod{p^2}, \label{wl0-4}\\
				&\sum_{k=1}^{\frac{p-1}{2}}\frac{1}{k^2}\equiv 0 \pmod{p}. \label{wl0-5}
			\end{align}		
	\end{lem}	
	\pf. Congruence  \eqref{wl0-3} is taken from \cite[(2.4)]{wt-6}, while \eqref{wl0-4} and \eqref{wl0-5} are from \cite{wt-7}.	\qed

\begin{lem} 
		Let $p>2$ be an odd prime. Then
		\begin{align}
				&\sum_{k=1}^{\frac{p-1}{2}}kH_{2k}
				\equiv  \frac{3p-2}{16}+\frac{2q_p(2)-pq_p(2)^2}{16}   \pmod{p^2}, \label{wl0-6}\\
				&\sum_{k=1}^{\frac{p-1}{2}}H_{2k}
				\equiv \frac{1-p}{2}-\frac{2q_p(2)-pq_p(2)^2}{4}    \pmod{p^2}, \label{wl0-7}
		\end{align}	
and
\begin{align*}
	\sum_{k=1}^{\frac{p-1}{2}}(2k+1)H_{2k}^2
	\equiv \frac{q_p(2)^2-2}{4}    \pmod{p}. 	
\end{align*}	
	\end{lem}	
	\pf. It is routine to see that
		\begin{align*}
			\begin{aligned}
		&\sum_{k=1}^{p-1}(-1)^kkH_k=\sum_{j=1}^{p-1}\sum_{k=j}^{p-1}(-1)^kk=\frac{(2p-1)H_{p-1}-H_{p-1}(-1)}{4}, \\
		&\sum_{k=1}^{p-1}kH_k=\frac{(2p^2-1)H_{p-1}-H_{p-1}(-1)-2+3p-p^2}{16}.
	\end{aligned}
\end{align*}	
By \eqref{wl0-1} and \eqref{wl0-4}, we have
\begin{align}
	\sum_{k=1}^{p-1}\frac{(-1)^k}{k}=\sum_{k=1}^{\frac{p-1}{2}}\frac{1}{k}-\sum_{k=1}^{p-1}\frac{1}{k}\equiv -2q_p(2)+pq_p(2)^2 \pmod{p^2}.  \label{wl0-8-0}
\end{align}			
Combining \eqref{wl0-1} and the above congruence yields	
\begin{align*}
	\sum_{k=1}^{\frac{p-1}{2}}kH_{2k}=\frac{\sum_{k=1}^{p-1}(1+(-1)^k)kH_k}{4}
	\equiv  \frac{3p-2}{16}+\frac{2q_p(2)-pq_p(2)^2}{16}   \pmod{p^2},		
\end{align*}		
as desired result. The proof of \eqref{wl0-7} is similar to that of \eqref{wl0-6}, the detailed steps are omitted.	

In view of  in \cite{wt-8}, we have
\begin{align}
\sum_{j=1}^{p-1}\sum_{i=1}^{j-1}\frac{(-1)^j}{ij}\equiv q_p(2)^2 \pmod{p}. \label{wl0-8-1}
\end{align}
By \eqref{wl0-2} and \eqref{wl0-4} we have
	\begin{align}
	H_{p-1}(-2)=\frac{1}{2}H_{\frac{p-1}{2}}(2)-H_{p-1}(2)\equiv 0 \pmod{p}. \label{wl0-8-2}
	\end{align}	
Using \eqref{wl0-1}, \eqref{wl0-2}, \eqref{wl0-4}, \eqref{wl0-8-1} and \eqref{wl0-8-2}, we obtain
\begin{align}
	\begin{aligned}
	\sum_{k=1}^{p-1}&(-1)^k(k+1)H_{k}^2=	\sum_{k=1}^{p-1}(-1)^k(k+1)\left(H_k(2)+2H(1,1;k)\right)\\
	&=2\sum_{j=1}^{p-1}\sum_{i=1}^{j}\frac{1}{ij}\sum_{k=j}^{p-1}(-1)^k(k+1)-\sum_{k=1}^{p-1}(-1)^k(k+1)H_k(2)\\
	&=\frac{(1+2p)}{4}\left(H_{p-1}^2+H_{p-1}(2)\right)+\frac{1}{2}\left(\sum_{j=1}^{p-1}\sum_{i=1}^{j-1}\frac{(-1)^j}{ij}+H_{p-1}(-2)\right)\\
	&+\frac{1}{2}H_{\frac{p-1}{2}}-\sum_{k=1}^{p-1}(-1)^k(k+1)H_k(2)\\
	&\equiv \frac{1}{2}q_p(2)^2-q_p(2)-\sum_{k=1}^{p-1}(-1)^k(k+1)H_k(2) \pmod{p}. \label{wl0-9}
\end{aligned}	
\end{align}			
Moreover, interchanging the order of summation and applying \eqref{wl0-2}, \eqref{wl0-8-0}, \eqref{wl0-8-2}  we get	
\begin{align}
	\sum_{k=1}^{p-1}(-1)^k(k+1)H_k(2)=\frac{(1+2p)H_{p-1}(2)+H_{p-1}(-2)+2H_{p-1}(-1)}{4}\equiv -q_p(2) \pmod{p}. \label{wl0-10}
\end{align}	
Addition, we also have	
\begin{align}
	\sum_{k=1}^{p-1}(k+1)H_k^2\equiv -\frac{1}{2}-\sum_{k=1}^{p-1}(k+1)H_k(2) \pmod{p}, \label{wl0-11}
\end{align}	
\begin{align}
	\sum_{k=1}^{p-1}(k+1)H_k(2)\equiv \frac{1}{2} \pmod{p}. \label{wl0-12}
\end{align}	
Since  \eqref{wl0-11} and \eqref{wl0-12} can be derived similarly to the proof of \eqref{wl0-9} and \eqref{wl0-10}, we omit their detailed proofs. Combining \eqref{wl0-9}--\eqref{wl0-12}, we deduce that	
\begin{align*}
	\sum_{k=1}^{\frac{p-1}{2}}(2k+1)H_{2k}^2&=\frac{1}{2}\sum_{k=1}^{p-1}(k+1)\left((-1^k+1)\right)H_k^2\\
	&\equiv \frac{q_p(2)^2-2}{4}    \pmod{p},
\end{align*}		
as claimed consequence. Thus, the proof of Lemma $2.2$ is finished.

\section{Proof of the Theorem 1.1}	
	According to the Dixon's theorem \cite{wt-2},
	\begin{align*}
		{}_3F_2\left[
		\begin{matrix}
			a,\ b,\ c \\
			1+a-b,\ 1+a-c
		\end{matrix}; 1
		\right]
		=
		\Gamma\left[
		\begin{matrix}
			1+\frac{a}{2},\ 1+a-b,\ 1+a-c,\ 1+\frac{a}{2}-b-c \\
			1+a,\ 1+\frac{a}{2}-b, \ 1+\frac{a}{2}-c ,\ 1+a-b-c
		\end{matrix}
		\right],
	\end{align*}	
	Feng and Hou \cite{wt-4} found a WZ pair $(F(n,k), G(n,k))$ as follows,
	\begin{align*}
		F(n,k)=\frac{(-1)^{n+k}n(2n+2k)!(2n-2k-2)!(4n)!(3n)!}{2^{12n}(n+k)!^2(2n-k-1)!(n-k-1)!^2n!^3(2n+k)!},
	\end{align*}	
	\begin{align*}
		G(n,k)=\frac{(-1)^{n+k}\alpha(n,k)(2n+2k)!(2n-2k)!(4n)!(3n)!}{2^{12n+8}(n+k)!^2(2n-k+1)!(n-k)!^2n!^3(2n+k+1)!},
	\end{align*}	
	where 
	\begin{align*}
		\begin{aligned}
			\alpha(n,k)&=16k^4-128k^2n^2+688n^4-76k^2n+988n^3-16k^2+508n^2+111n+9.
		\end{aligned}
	\end{align*}
	which satisfy the WZ equation
	\begin{align}
		F(n+1,k)-F(n,k)=G(n,k+1)-G(n,k).  \label{3pr-1}
	\end{align}
	
	Our proof will start from equation \eqref{3pr-1}. Summing up \eqref{3pr-1} for $k$ from $0$ to $p-1$, we derive that
	\begin{align}
		\sum_{k=0}^{p-1}F(n+1,k)-\sum_{k=0}^{p-1}F(n,k)=G(n,p)-G(n,0).  \label{3pr-2}
	\end{align}
	Since $G(n,k)=0$ for $n<k$ and $F(0,k)=0$, we continue to sum both sides of \eqref{3pr-2} over $n$ from $0$ to $\frac{p-3}{2}$, which yields	
	\begin{align*}
		\sum_{k=0}^{p-1}F(\frac{p-1}{2},k)=-\sum_{n=0}^{\frac{p-3}{2}}G(n,0).
	\end{align*}
	
	Note that	
	\begin{align*}
		&G(n,0)=\frac{(-1)^n(172n^2+75n+9)}{2^{12n+8}}{2n \choose n}^3{3n \choose n}{4n\choose 2n},
	\end{align*}
	and	
	\begin{align*}
		\begin{aligned}
		F(\frac{p-1}{2},k)&=(-1)^{\frac{p-1}{2}}\frac{p-1}{2^{6p-7}}{ \frac{3p-3}{2} \choose  \frac{p-1}{2}}{p-1\choose \frac{p-1}{2}}\\
		&\times(-1)^k(p-1-k){p-1+2k \choose \frac{p-1+2k}{2}}{p-3-2k\choose \frac{p-3-2k}{2}}{2p-2\choose p-1+k}.
	\end{aligned}
	\end{align*}	
	Consequently, with some necessary arrangements, we obtain	
	\begin{align}
	\begin{aligned}
		\sum_{k=0}^{\frac{p-3}{2}}&\frac{(-1)^k(172k^2+75k+9)}{2^{12k}}{2k \choose k}^3{3k \choose k}{4k\choose 2k}\\
		&=-(-1)^{\frac{p-1}{2}}\frac{2^5p(p-1)^3}{2^{6(p-1)}(2p-1)(p-2)}{ \frac{3p-3}{2} \choose  \frac{p-1}{2}}{p-1\choose \frac{p-1}{2}}^3{2p-1\choose p-1}\\
		&-(-1)^{\frac{p-1}{2}}\frac{2^5(p-1)}{2^{6(p-1)}}{ \frac{3p-3}{2} \choose  \frac{p-1}{2}}{p-1\choose \frac{p-1}{2}}\mathcal{M}, \label{3pr-3}
	\end{aligned}
\end{align}
where 
	\begin{align*}
	\begin{aligned}
	\mathcal{M}=\sum_{k=1}^{\frac{p-3}{2}}(-1)^k\frac{(p-1-k)(p-1-2k)}{p-2-2k}{p-1+2k \choose \frac{p-1+2k}{2}}{p-1-2k\choose \frac{p-1-2k}{2}}{2p-2\choose p-1+k}. 
	\end{aligned}
\end{align*}	
Next, we first evaluate $\mathcal{M}$ module $p^4$.
		
For a nonnegative integer $k$, we have
	\begin{align}
	\begin{aligned}
	{p-1+2k\choose \frac{p-1+2k}{2}}&=p{p-1 \choose \frac{p-1}{2}}\frac{2^{2k}\prod_{j=1}^{2k-1}(p-1+j)}{\prod_{j=1}^{k}(p+2j-1)^2}\\
	&\equiv {p-1 \choose \frac{p-1}{2}}\frac{4^{2k}p}{(p+2k){2k \choose k}}\frac{1+\frac{p}{2}H_k+\frac{p^2}{4}H(1,1;k)}{1+pH_{2k}+p^2H(1,1,;2k)}\\
	&\equiv \left(1+p(H_k-H_{2k})+\frac{p^2}{4}(2H_k^2+2H_{2k}^2+2H_{2k}(2)^2-H_k(2)-H_kH_{2k})\right)\\
	&\times {p-1 \choose \frac{p-1}{2}}\frac{4^{2k}p}{(p+2k){2k \choose k}} \pmod {p^4}, 
	\label{3pr-4}
	\end{aligned}
\end{align}
and
	\begin{align}
	\begin{aligned}
		{p-1-2k\choose \frac{p-1-2k}{2}}&= \frac{{p-1\choose (p-1)/2}{(p-1)/2 \choose k}^2}{{2k\choose k}{p-1\choose 2k}}\\
		&\equiv {p-1\choose \frac{p-1}{2}}\frac{{2k\choose k}}{16^k}\frac{(1-pH_{2k}+p^2H(1,1;2k))}{(1-\frac{p}{2}H_k+\frac{p^2}{4}H(1,1;k))}\\
		&\equiv \left(1+p(H_k-H_{2k})+\frac{p^2}{4}(2H_k^2+2H_{2k}^2+H_k(2)-2H_{2k}(2)-H_kH_{2k})\right)\\
		&\times 
		{p-1\choose \frac{p-1}{2}}\frac{{2k\choose k}}{16^k}
		 \pmod {p^3}. \label{3pr-5}
	\end{aligned}
\end{align}
Meanwhile, since
	\begin{align}
	\begin{aligned}
	{p-1\choose k}=\frac{(p-1)(p-2)\cdots(p-k)}{k!}\equiv (-1)^k\left(1-pH_k+p^2H(1,1;k)\right)  \pmod {p^3}. \label{3pr-5-0}
	\end{aligned}
\end{align}
Therefore, by \eqref{wl0-1}, \eqref{wl0-2} and \eqref{3pr-5-0} yields
	\begin{align}
	\begin{aligned}
		&{2p-2\choose p-1-k}=\frac{(p+p-2)(p+p-3)\cdots(p+k)}{(p-1-k)!}\\
		&\equiv \frac{k(1+p\sum_{k\leq j \leq p-2}\frac{1}{j}+p^2\sum_{k\leq i<j\leq p-2}\frac{1}{ij})}{p-1}{p-1\choose k}\\
		&\equiv -k(-1)^k\left(1+2p+4p^2+\frac{p+2p^2}{k}-(2p+4p^2+\frac{2p^2}{k})H_k+2p^2H_k^2\right)  \pmod {p^3}. \label{3pr-6}
	\end{aligned}
\end{align}

Moreover, we deduce that for   $1\leq k\leq (p-3)/2$,  modulo $p^3$
\begin{align}
	\begin{aligned}
	\frac{k(p-1-2k)(p-1-k)}{(p-2-2k)(p+2k)}\equiv -\left(\frac{2k+1-3p}{4}-\frac{2p-5p^2}{16k}+\frac{p^2}{16k^2}+\frac{2p+p^2}{16(k+1)}+\frac{p^2}{16(1+k)^2}\right). \label{3pr-7}
	\end{aligned}
\end{align}
Note that
\begin{align}
	\begin{aligned}
	\sum_{k=1}^{\frac{p-1}{2}}\frac{1}{2k-1}=\sum_{k=1}^{\frac{p-1}{2}}\frac{1}{p-2k}\equiv -\frac{1}{2}H_{\frac{p-1}{2}}	\pmod{p}. \label{3pr-8}
	\end{aligned}
\end{align}

Combining \eqref{wl0-1}, \eqref{3pr-4}--\eqref{3pr-8}, and after explicit calculation, we obtain
\begin{align*}
	\begin{aligned}
		\mathcal{M}&/{p-1\choose \frac{p-1}{2}}^2\equiv -\frac{3p^3+3p}{16}+\frac{2p^2-3p^3}{8}H_{\frac{p-1}{2}}+\frac{p^3}{8}H_{\frac{p-1}{2}}(2)-\frac{p^3}{2}\sum_{k=1}^{\frac{p-1}{2}}\frac{H_{2k}}{k}\\
		&-(p^2-2p^3)\sum_{k=1}^{\frac{p-1}{2}}kH_{2k}-\frac{p^2-p^3}{2}\sum_{k=1}^{\frac{p-1}{2}}H_{2k}+\frac{p^3}{2}(2k+1)H_{2k}^2 \pmod{p^4}. 
	\end{aligned}
\end{align*}
It is well-know the congruence of Morley \cite{wt-9}: for prime $p>2$,
\begin{align}
{p-1\choose \frac{p-1}{2}}\equiv (-1)^{\frac{p-1}{2}}4^{p-1} \pmod{p^3}. \label{3pr-10}
\end{align}
Using Lemma $1.1$, Lemma $1.2$, \eqref{3pr-10} and  simplifying, we arrive at
\begin{align}
	\begin{aligned}
		\mathcal{M}\equiv -\frac{3p+2p^2+6p^3}{16}-\frac{9p^2-2p^3}{8}q_p(2)-\frac{45p^3}{16}q_p(2)^2\pmod{p^4}. \label{3pr-11}
	\end{aligned}
\end{align}

By \eqref{wl0-4} and \eqref{wl0-5}, we have
\begin{align}
	\begin{aligned}
	{\frac{3p-3}{2}\choose \frac{p-1}{2}}&=\frac{p(p+\frac{p-3}{2})(p+\frac{p-5}{2})\cdots(p+1)}{(p-1)!}\\
	&\equiv \frac{2p}{3p-1}\left(1+pH_{\frac{p-1}{2}}+\frac{p^2}{2}H(1,1;\frac{p-1}{2})\right)\\
&\equiv -2p-6p^2-18p^3+4(p^2+3p^3)q_p(2)-6p^3q_p(2)^2 \pmod{p^4}, \label{3pr-12}
	\end{aligned}
\end{align}
and 
\begin{align}
{2p-1\choose p-1}\equiv 1\pmod{p^3}. \label{3pr-13}
\end{align}
For positive integer $\alpha$, we have
\begin{align*}
	(2^{p-1})^\alpha=(1+pq_p(2))^\alpha.
\end{align*}

Then substituting \eqref{3pr-10}--\eqref{3pr-13} into \eqref{3pr-3} and simplifying, we get
\begin{align*}
	\begin{aligned}
		\sum_{k=0}^{\frac{p-3}{2}}&\frac{(-1)^k(172k^2+75k+9)}{2^{12k}}{2k \choose k}^3{3k \choose k}{4k\choose 2k}\\
		&\equiv -20p^2-48p^3-152p^4+(64p^3+96p^4)q_p(2)-96p^4q_p(2)^2 \pmod{p^5}.
	\end{aligned}
\end{align*}
With the help of \eqref{3pr-10}, \eqref{3pr-12} and \eqref{3pr-13} gives
\begin{align*}
	\begin{aligned}
	(-1)^{\frac{p-1}{2}}\frac{43p(9+75(p-1)/2+(p-1)^2)}{2^{6(p-1)}(2p-1)}{p-1\choose (p-1)/2}{(3p-3)/2\choose (p-1)/2}{2p-1\choose p-1}\\
	\equiv 29p^2+48p^3+152p^4-(58p^3+96p^4)q_p(2)+87p^4q_p(2)^2  \pmod{p^5}.
	\end{aligned}
\end{align*}
It follows that
\begin{align*}
	\begin{aligned}
		\sum_{k=0}^{\frac{p-1}{2}}&\frac{(-1)^k(172k^2+75k+9)}{2^{12k}}{2k \choose k}^3{3k \choose k}{4k\choose 2k}\\
		&=(-1)^{\frac{p-1}{2}}\frac{43p(9+75(p-1)/2+(p-1)^2)}{2^{6(p-1)}(2p-1)}{p-1\choose (p-1)/2}{(3p-3)/2\choose (p-1)/2}{2p-1\choose p-1}\\
		&+\sum_{k=0}^{\frac{p-3}{2}}\frac{(-1)^k(172k^2+75k+9)}{2^{12k}}{2k \choose k}^3{3k \choose k}{4k\choose 2k}\\
		&\equiv 9p^2+6p^3q_p(2)-9p^4q_p(2)^2 \pmod{p^5}.
	\end{aligned}
\end{align*}
This concludes the proof. \qed

	\section{Proof of the Theorem 1.2}	
	
	Summing both sides of \eqref{3pr-1} over $n$ from $0$ to $p-1$ yields
	\begin{align*}
		F(p,k)-F(0,k)=\sum_{n=0}^{p-1}G(p,k+1)-\sum_{n=0}^{p-1}G(p,k).
	\end{align*}	
	Continuing to sum over $k$ from $0$ to $p-1$, and utilizing  $G(n,k)=0$ for $n<k$ together with $F(0,k)=0$, we obtain	
	\begin{align*}
		\sum_{k=0}^{p-1}F(p,k)=-\sum_{n=0}^{p-1}G(n,0),
	\end{align*}		
where
\begin{align*}
	F(p,k)=-\frac{(-1)^k6kp(2p-k)}{2^{12p}}{3p\choose p}{2p\choose p}{2p+2k\choose p+k}{2p-2k-2\choose p-k-1}{4p\choose 2p+k}.
\end{align*}
After some necessary arrangement, we have
	\begin{align}
		\begin{aligned}
	\sum_{k=0}^{p-1}&\frac{(-1)^k(172k^2+75k+9)}{2^{12k}}{2k \choose k}^3{3k \choose k}{4k\choose 2k}\\
	&=\frac{2^8p}{2^{12p}}{3p\choose p}{2p\choose p}\sum_{k=0}^{p-1}(-1)^k(2p-k){2p+2k\choose p+k}{2p-2k-2\choose p-k-1}{4p\choose 2p+k}. \label{4pr-1}
\end{aligned}	
\end{align}
Now we split the sum on the right-hand of \eqref{4pr-1} into five pieces:	
	\begin{align}
	\begin{aligned}
		\sum_{k=0}^{p-1}&(-1)^k(2p-k){2p+2k\choose p+k}{2p-2k-2\choose p-k-1}{4p\choose 2p+k}\\
		&=\sum_{k=1}^{\frac{p-3}{2}}G^*(k)+ \sum_{k=\frac{p+1}{2}}^{p-2}G^*(k)+G^*(0)+G^*(\frac{p-1}{2})+G^*(p-1)\label{4pr-2},
	\end{aligned}	
\end{align}
where	
\begin{align*}
G^*(k)=(-1)^k\frac{(2p-k)(p-k)}{2(2p-2k-1)}{2p+2k\choose p+k}{2p-2k\choose p-k}{4p\choose 2p+k}.
\end{align*}	

On the one hand, for $1\leq k \leq (p-3)/2$, we have	
\begin{align*}
	{4p\choose 2p+k}=\frac{4p(3p+p-1)\cdots(3p+1)3p(2p+p-1)\cdots(2p+k+1)}{(p+p-k)!}\\
	\equiv \frac{12p}{p-k}\frac{(1+3pH_{p-1})(1+2p\sum_{i=k+1}^{p-1}\frac{1}{i})}{1+pH_{p-k}}{p-1\choose k} \pmod{p^3}.
\end{align*}	
In a similar way, we obtain
\begin{align}
	{2p-2k\choose p-k}\equiv  -\frac{2p}{k{2k\choose k}}\left(1-2pH_{k-1}+2pH_{2k-1}\right)   \pmod{p^3}, \label{4pr-2-1}
\end{align}
	\begin{align*}
	{2p+2k\choose p+k}\equiv \frac{(-1)^k2p}{(2k+1){p+k\choose 2k+1}}\left(1+2pH_{2k}-2pH_{p-k-1}\right) \pmod{p^2},
\end{align*}
and
\begin{align*}
	{p+k\choose 2k+1}\equiv \frac{(-1)^kp}{(2k+1){2k\choose k}} \pmod{p^3}.
\end{align*}
these three congruences,  together with \eqref{wl0-1} and \eqref{3pr-5-0} gives 
\begin{align}
	\begin{aligned}
	(-1^k)&\frac{p-k}{2}{2p+2k\choose p+k}{2p-2k\choose p-k}{4p\choose 2p+k}\\
	&\equiv -\frac{24p^2}{k}\left(1-2pH_{p-1-k}-pH_{p-k}-5pH_k+4pH_{2k-1}+\frac{3p}{k}\right) \pmod{p^4}. \label{4pr-3}
\end{aligned}
\end{align}
It is easy for us to check that
\begin{align}
	\frac{2p-k}{2p-2k-1}\equiv \frac{2kp}{(2k+1)^2}-\frac{2p-k}{2k+1}  \pmod{p^2}.    \label{4pr-4}
\end{align}
Thus, combining \eqref{4pr-3}, \eqref{4pr-4} and the fact $H_{p-1-k}\equiv H_k \quad \pmod{p}$, we have
\begin{align}
	\begin{aligned}
   \sum_{k=1}^{\frac{p-3}{2}}G^*(k)&\equiv  -48p^3H_{\frac{p-3}{2}}+(96p^3-24p^2)\sum_{k=1}^{\frac{p-3}{2}} \frac{1}{2k+1}-48p^3 \sum_{k=1}^{\frac{p-3}{2}} \frac{1}{(2k+1)^2} \\
   &+192p^3 \sum_{k=1}^{\frac{p-3}{2}} \frac{H_k}{2k+1}-96p^3\sum_{k=1}^{\frac{p-3}{2}} \frac{H_{2k-1}}{2k+1} \pmod{p^4}.   \label{4pr-5}
\end{aligned}
\end{align}

On the other hand, substituting $k$ with $p-1-k$, we get
\begin{align*}
	\begin{aligned}
		\sum_{k=\frac{p+1}{2}}^{p-2}G(k)^*&=\sum_{k=1}^{\frac{p-3}{2}}\frac{(-1)^k(p+1+k)(2p-k)}{2(4p-2k-1)}{4p-2k \choose 2p-k}{2k\choose k}{4p\choose p+1+k}.
	\end{aligned}
\end{align*}
Meanwhile, for $1\leq k \leq \frac{p-3}{2}$ we have
\begin{align*}
	\begin{aligned}
	{4p\choose p+1+k}&=\frac{4p(3p+p+1)\cdots(3p+1)3p(3p-1)\cdots(3p-k)}{(p+1+k)!}\\
	&\equiv \frac{(-1)^k12p(1+3pH_{p-1})(1-2pH_k)}{(k+1)(1+pH_{k+1})} \pmod{p^3}.
	\end{aligned}
\end{align*}
and
\begin{align*}
	\begin{aligned}
		{4p-2k \choose 2p-k}{2k\choose k}&={4p\choose 2p}{2p\choose k}^2/{4p\choose 2k}
		\equiv -{4p\choose 2p}\frac{2p(1+p\sum_{i=p+1-k}^{p-1}\frac{1}{i})^2(1-pH_{k-1})^2}{k(1+3p\sum_{i=p+1-2k}^{p-1}\frac{1}{i})(1-pH_{2k-1})} \pmod{p^3}.
	\end{aligned}
\end{align*}
After careful simplification and  combination, we obtain
\begin{align}
	\begin{aligned}
		(-1)^k&{4p\choose p+1+k}{4p-2k \choose 2p-k}{2k\choose k}/{4p\choose 2p}\\
		&\equiv -24p^2\sum_{k=1}^{\frac{p-3}{2}}\frac{1}{k(k+1)}\left(1+4pH_{2k-1}-8pH_{k-1}-\frac{4p}{k}-\frac{p}{k+1}\right) \pmod{p^4}.   \label{4pr-6}
	\end{aligned}
\end{align}
Obviously, from \eqref{wl0-1} we have
\begin{align}
	\begin{aligned}
	{4p \choose 2p}=\frac{4p(3p+p-1)\cdots(3p)(2p+p-1)\cdots(2p+1)}{(2p)!}\equiv 6    \pmod{p^3},   \label{4pr-7}
	\end{aligned}
\end{align}
and
\begin{align}
	\begin{aligned}
	\frac{(p+1+k)(2p-k)}{4p-2k-1}
		\equiv -p\left(\frac{3}{2(2k+1)}+\frac{1}{(2k+1)^2}-\frac{1}{2}\right)+\frac{2k+1}{4}-\frac{1}{4(2k+1)}     \pmod{p^2}.   \label{4pr-8}
	\end{aligned}
\end{align}
Combining \eqref{4pr-6}--\eqref{4pr-8}, we arrive at
\begin{align}
	\begin{aligned}
	\sum_{k=\frac{p+1}{2}}^{p-2}G^*(k)&\equiv 72(4p^3-p^2)\sum_{k=1}^{\frac{p-3}{2}}\frac{1}{1+2k}
	-288p^3\sum_{k=1}^{\frac{p-3}{2}}\frac{1}{(1+2k)^2}
	-144p^3\sum_{k=1}^{\frac{p-3}{2}}\frac{1}{k}\\
	&-288p^3\sum_{k=1}^{\frac{p-3}{2}}\frac{H_{2k-1}}{2k+1}+576p^3\sum_{k=1}^{\frac{p-3}{2}}\frac{H_k}{2k+1}  \pmod{p^4}. \label{4pr-9}
	\end{aligned}
\end{align}

Furthermore, by \eqref{wl0-4}, \eqref{wl0-5} we have
\begin{align}
	\begin{aligned}
	\sum_{k=1}^{\frac{p-3}{2}}\frac{1}{1+2k}=\sum_{k=1}^{\frac{p-1}{2}}\frac{1}{p-2k}-1\equiv q_p(2)-\frac{p}{2}q_p(2)^2-1 \pmod{p^2},  \label{4pr-9-1}
	\end{aligned}
\end{align}
and
\begin{align}
	\begin{aligned}
		\sum_{k=1}^{\frac{p-3}{2}}\frac{1}{(1+2k)^2}=\sum_{k=1}^{\frac{p-1}{2}}\frac{1}{(p-2k)^2}-1\equiv -1 \pmod{p}.  \label{4pr-9-2}
	\end{aligned}
\end{align}
In view of \cite[(2.1)]{wt-6} and \cite[(3.11)]{wt-10}: for prime $p$,
\begin{align*}
H_{\frac{p-1}{2}-k}\equiv H_{\frac{p-1}{2}}+2H_{2k}-H_k \pmod{p},
\end{align*}
and
\begin{align*}
	\sum_{k=1}^{\frac{p-1}{2}}\frac{H_k}{k}\equiv 2q_p(2)^2 \pmod{p}.
\end{align*}
Using \eqref{wl0-3}, \eqref{wl0-4}, \eqref{4pr-9-1}  and the above congruences, we obtain
\begin{align}
	\sum_{k=1}^{\frac{p-3}{2}}\frac{H_k}{2k+1}=\sum_{k=1}^{\frac{p-1}{2}}\frac{H_{\frac{p-1}{2}-k}}{p-2k}\equiv -2q_p(2)^2 \pmod{p}, \label{4pr-9-3}
\end{align}
and
\begin{align}
	\sum_{k=1}^{\frac{p-3}{2}}\frac{H_{2k-1}}{2k+1}&=\sum_{k=1}^{\frac{p-1}{2}}\frac{H_{p-1-2k}}{p-2k}+\sum_{k=1}^{\frac{p-3}{2}}\frac{1}{2k+1}-\frac{1}{2}H_{\frac{p-1}{2}}+\frac{1}{p-1}    \\
	&\equiv 2q_p(2)-\frac{1}{2}q_p(2)^2-2 \pmod{p}. \label{4pr-9-4}
\end{align}

In addition, with the help of \eqref{wl0-1}, \eqref{wl0-2}, \eqref{wl0-4} and \eqref{wl0-5} yields
\begin{align}
	\begin{aligned}
	{4p\choose \frac{5p-1}{2}}
	&=\frac{4p(3p+p-1)\cdots(3p)(2p+p-1)\cdots(2p+\frac{p+1}{2})}{(p+\frac{p+1}{2})\cdots(p+1)p!}\\
	&\equiv {p-1\choose \frac{p-1}{2}}\left(24p-72p^2+216p^3+(144p^2-432p^3)q_p(2)+360p^3q_p(2)^2\right) \pmod{p^4}.\label{4pr-10}
\end{aligned}
\end{align}
Similarly, we also deduce that
\begin{align}
	\begin{aligned}	
		{3p-1\choose \frac{3p-1}{2}}\equiv 2{p-1\choose \frac{p-1}{2}}\left(1+4pq_p(2)+6p^2q_p(2)^2\right) \pmod{p^3},\label{4pr-11}
	\end{aligned}	
\end{align}
and
\begin{align}
	{4p\choose p}\equiv  2{2p\choose p}\equiv 4{2p-1\choose p-1} \equiv4 \pmod{p^3}.\label{4pr-12}
\end{align}
Then, utilizing \eqref{3pr-10}, \eqref{4pr-7}, \eqref{4pr-10}--\eqref{4pr-12} and simplifying we get
\begin{align*}
	G^*(0)\equiv -24p^2\left(1+2p\right) \pmod{p^4},
\end{align*}
\begin{align*}
	G^*(p-1)\equiv -72p^2\left(1+4p\right) \pmod{p^4},
\end{align*}
and
\begin{align*}
	G^*(\frac{p-1}{2})\equiv 2^{6(p-1)}\left(24p+240p^2q_p(2)+1080p^3q_p(2)^2\right) \pmod{p^4}.
\end{align*}
The above three congruences, together with \eqref{wl0-4}, \eqref{4pr-5}, \eqref{4pr-9}, \eqref{4pr-9-1}--\eqref{4pr-9-4} gives
\begin{align}
	\begin{aligned}
		\sum_{k=0}^{p-1}G^*(k)\equiv 24p+288p^2q_p(2)+1584p^3q_p(2)^2  \pmod{p^4}. \label{4pr-18}
	\end{aligned}
\end{align}
It is clear that
\begin{align}
	{3p\choose p}{2p \choose p}\equiv 6 \pmod{p^3}. \label{4pr-19}
\end{align}
 Finally, substituting \eqref{4pr-18} and \eqref{4pr-19} into \eqref{4pr-1} and simplifying we arrive at
\begin{align*}
	\sum_{k=0}^{p-1}&\frac{(-1)^k(172k^2+75k+9)}{2^{12k}}{2k \choose k}^3{3k \choose k}{4k\choose 2k}\equiv 9p^2 \pmod{p^5},
\end{align*}
as claimed result. Thus, we complete the proof of Theorem $1.2$. \qed

	\section{Proof of the Theorem 1.3}
Applying the following ${}_7F_6$-series identities \cite{wt-11}:
\begin{align*}
	\begin{aligned}
	{}_7F_6&\left[
	\begin{matrix}
		a-\frac{1}{2},\ \frac{2a+2}{3},\ 2b-1,\ 2c-1,\ 2+2a-2b-2c,\ a+n,\ -n \\
		\frac{2a-1}{3},\ 1+a-b,\ 1+a-c.\ b+c-\frac{1}{2},\ 2a+2n.\ -2n
	\end{matrix}; 1
	\right]\\
	&=\frac{\left(\frac{1}{2}+a\right)_n\left(b\right)_n\left(c\right)_n\left(a+b-c+\frac{3}{2}\right)_n}{\left(\frac{1}{2}\right)_n\left(1+a-b\right)_n\left(1+a-c\right)_n\left(b+c-\frac{1}{2}\right)_n}.
\end{aligned}
\end{align*}
Feng and Hou \cite{wt-4} found a WZ pair $(F^{'}(n,k), G^{'}(n,k))$ as follows,
\begin{align*}
	F^{'}(n,k)=\frac{(-1)^{k}n(3n+3k-1)(2n+2k-2)!(3n+k-1)!(2n)!k!}{2^{6n+2k}(n+k-1)(n+k)!^2(n-k-1)!(2n+k-1)!n!^2(2k)!},
\end{align*}	
\begin{align*}
	G^{'}(n,k)=\frac{(-1)^{k+1}\beta(n,k)(2n+2k-3)!(3n+k-1)!(2n)!(k-1)!}{2^{6n+2k+2}(n+k)!^2(n-k)!(2n+k)!n!^2(2k-2)!},
\end{align*}	
where 
\begin{align*}
	\begin{aligned}
		\beta(n,k)&=22n^3+(32k-3)n^2+(10n^2+2n-3)n+k^2-k.
	\end{aligned}
\end{align*}
which satisfy the WZ equation
\begin{align}
	F^{'}(n+1,k)-F^{'}(n,k)=G^{'}(n,k+1)-G^{'}(n,k).  \label{5pr-1}
\end{align}

Summing both sides of \eqref{5pr-1} over $k$ from $0$ to $p-1$ gives 
\begin{align*}
	\sum_{k=0}^{p-1}F^{'}(n+1,k)-\sum_{k=0}^{p-1}F^{'}(n,k)=G^{'}(n,p)-G^{'}(n,0). 
\end{align*}
Notice that  for $n<k$,  $G^{'}(n,k)=0$ and $F^{'}(0,k)=0$. we   proceed  to sum over $n$ from $2$ to $p-1$, 
\begin{align*}
	\sum_{k=0}^{p-1}F^{'}(p,k)-\sum_{k=0}^{p-1}F^{'}(2,k)=-\sum_{n=2}^{p-1}G^{'}(n,0), 
\end{align*}
after some necessary arrangement and calculation, this identity is equivalent to
\begin{align}
	\begin{aligned}
	\frac{p^2}{2^{6p+1}}{2p\choose p}&\sum_{k=0}^{p-1}\frac{(-1)^k(3p+3k-1)}{2^{2k}(p+k-1)(p+k)(2p+2k-1)}\frac{{2p+2k\choose p+k}{3p+k-1\choose p}{p-1\choose k}}{{2k\choose k}}\\
	&=\frac{15}{128}-\frac{1}{24}\sum_{k=2}^{p-1}\frac{22k^2-3k-3}{2^{6k}k(k+1)(2k-1)}{2k\choose k}^2{3k\choose k}.  \label{5pr-2}
\end{aligned}
\end{align}
Let 
\begin{align*}
	\begin{aligned}
	g_n=\frac{11n+3}{2^{6n}}{2n\choose n}^2{3n\choose n},
	\end{aligned}
\end{align*}
which satisfies
\begin{align}
	\begin{aligned}
		g_n-48G^{'}(n,0)=\Delta_n\left(\frac{n^3(2n)!(3n)!}{2^{6n-5}(2n-1)(n-1)n!^5}\right), \label{5pr-2-1}
	\end{aligned}
\end{align}
where $\Delta_n$ denote the forward difference operator on $n$: $\Delta_nf(n)=f(n+1)-f(n)$. Then summing both sides of \eqref{5pr-2-1} over $n$ from $2$ to $p-1$, we obtain
\begin{align}
	\begin{aligned}
		\sum_{n=2}^{p-1}&\frac{11n+3}{2^{6n}}{2n\choose n}^2{3n\choose n}
		-48\sum_{n=2}^{p-1}\frac{22n^2-3n-3}{2^{6n+3}3n(n-1)(2n-1)}{2n\choose n}^2{3n\choose n}\\
		&=\frac{p^3}{2^{6p-5}(p-1)(2p-1)}{2p\choose p}^2{3p\choose p}-\frac{45}{4}.  \label{5pr-3}
	\end{aligned}
\end{align}

Next, we evaluate sum the left-hand of \eqref{5pr-2} modulo $p$. Firstly, we perform the replacement $k\rightarrow p-k$,
\begin{align}
	\begin{aligned}
		&\sum_{k=1}^{p-1}\frac{(-1)^k(3p+3k-1)}{2^{2k}(p+k-1)(p+k)(2p+2k-1)}\frac{{2p+2k\choose p+k}{3p+k-1\choose p}{p-1\choose k}}{{2k\choose k}}\\
		&=-\frac{1}{4^p}\sum_{k=1}^{p-1}\frac{(-1)^k2^{2k}(6p-3k-1)k}{(2p-k-1)(2p-k)(4p-2k-1)(p-k)}\frac{{4p-2k\choose 2p-k}{4p-k-1\choose p}{p-1\choose k}}{{2p-2k\choose p-k}}. \label{5pr-4}
	\end{aligned}
\end{align}
We split sums  the right-hand of \eqref{5pr-4} into four pieces
\begin{align}
	\begin{aligned}
	\sum_{k=1}^{p-1}&\frac{(-1)^k2^{2k}(6p-3k-1)k}{(2p-k-1)(2p-k)(4p-2k-1)(p-k)}\frac{{4p-2k\choose 2p-k}{4p-k-1\choose p}{p-1\choose k}}{{2p-2k\choose p-k}}\\
	&= \sum_{k=1}^{\frac{p-3}{2}}F^*(p,k)+\sum_{k=\frac{p+1}{2}}^{p-2}F^*(p,k)+F^*(p,p-1)+F^*(p,\frac{p-1}{2}), \label{5pr-5}
	\end{aligned}
\end{align}
where 
\begin{align*}
	\begin{aligned}
		F^*(p,k)=\frac{(-1)^k2^{2k}(6p-3k-1)k}{(2p-k-1)(2p-k)(4p-2k-1)(p-k)}\frac{{4p-2k\choose 2p-k}{4p-k-1\choose p}{p-1\choose k}}{{2p-2k\choose p-k}}.
	\end{aligned}
\end{align*}

For $1\leq k \leq (p-3)/2$, we have
\begin{align*}
	\begin{aligned}
	{4p-k-1\choose p}=\frac{(3p+p-1)3p(3p-1)\cdots(3p-k)}{p!}
	\equiv 3(-1)^k/{p-1\choose k} \pmod{p},
	\end{aligned}
\end{align*}
It is easy for us to see that
\begin{align}
	\begin{aligned}
		{4p\choose 2p}/{2p\choose p}=3 \pmod{p}. \label{5pr-6-0}
	\end{aligned}
\end{align}
For positive integer $n$, we have
\begin{align}
	\begin{aligned}
{2n-2k\choose n-k}=\frac{{2n\choose n}{n\choose k}^2}{{2k\choose k}{2n\choose 2k}}.	\label{5pr-6}	
	\end{aligned}
\end{align}
Applying \eqref{5pr-6-0} and \eqref{5pr-6}, we obtain
\begin{align*}
	\begin{aligned}
		{4p-2k\choose 2p-k}/{2p-2k\choose p-k}=\frac{{4p \choose 2p}{2p \choose 2k}{2p \choose k}^2}{{2p \choose p}{4p \choose 2k}{p \choose k}^2}\equiv 6 \pmod{p}.
	\end{aligned}
\end{align*}
It follows that
\begin{align}
	\begin{aligned}
	 \sum_{k=1}^{\frac{p-3}{2}}F^*(p,k)
	 &\equiv 18\sum_{k=1}^{\frac{p-3}{2}} \frac{2^{2k}(6p-3k-1)k}{(2p-k-1)(2p-k)(4p-2k-1)(p-k)}\\
	 &\equiv -18\sum_{k=1}^{\frac{p-3}{2}} 4^k\left(\frac{1}{k}+\frac{2}{2k+1}+\frac{2}{k+1}\right)\\
	 &\equiv -36-9\sum_{k=1}^{\frac{p-1}{2}}\frac{4^k}{2k-1}-9\sum_{k=1}^{\frac{p-1}{2}} \frac{4^k}{k}  \pmod{p}. \label{5pr-8}
	\end{aligned}
\end{align}
Moreover, it is clear for us to check that
\begin{align*}
	\begin{aligned}
		{2p\choose 2k+2}/{4p\choose 2(p-1-k)}\equiv \frac{1}{6} \pmod{p},
	\end{aligned}
\end{align*}
\begin{align*}
	\begin{aligned}
		{2p\choose p-1-k}^2/{p\choose k+1}^2\equiv 4 \pmod{p},
	\end{aligned}
\end{align*}
and
\begin{align*}
	\begin{aligned}
		{3p+k\choose p}{p-1\choose k}\equiv 3(-1)^k \pmod{p}.
	\end{aligned}
\end{align*}
Combining \eqref{5pr-6-0}, \eqref{5pr-6}, and the above three congruence, we get
\begin{align}
	\begin{aligned}
		\sum_{k=\frac{p+1}{2}}^{p-2}F^*(p,k)&=4^{p-1}\frac{{4p\choose 2p}}{{2p\choose p}}\sum_{k=1}^{\frac{p-3}{2}}\frac{(-1)^k(3p+3k+2)(p-1-k)}{4^k(p+k)(p+k+1)(2p+2k+1)(1+k)}\\
		&\times\frac{{2p\choose 2k+2}{2p\choose p-1-k}^2{3p+k\choose p}{p-1\choose k}}{{4p\choose 2(p-1-k)}{p\choose k+1}^2}\\
		&\equiv -6\sum_{k=1}^{\frac{p-3}{2}}\frac{1}{4^k}\left(\frac{2}{k}-\frac{2}{2k+1}-\frac{1}{k+1}\right)	 \\
		&\equiv 12\sum_{k=1}^{\frac{p-1}{2}}\frac{1}{4^kk}+48\sum_{k=1}^{\frac{p-1}{2}}\frac{1}{4^k(2k-1)}-42 \pmod{p},
	 \label{5pr-9}
	\end{aligned}
\end{align}
where we used \eqref{5pr-6} and  made the replacement $k\to p-k$ in the first step.

Meanwhile, applying Fermat's Little Theorem we have
\begin{align}
	\begin{aligned}
	\sum_{k=1}^{\frac{p-1}{2}}\frac{1}{4^kk}=\frac{1}{2^{p}}\sum_{k=1}^{\frac{p-1}{2}}
	\frac{4^k}{p+1-2k}\equiv 	-\frac{1}{2}\sum_{k=1}^{\frac{p-1}{2}}\frac{1}{4^k(2k-1)} \pmod{p},  \label{5pr-10}
	\end{aligned}
\end{align}
and
\begin{align}
	\begin{aligned}
		\sum_{k=1}^{\frac{p-1}{2}}\frac{1}{(2k-1)4^k}=\frac{1}{2^{p+1}}\sum_{k=1}^{\frac{p-1}{2}}\frac{4^k}{p-2k} \equiv 	-\frac{1}{8}\sum_{k=1}^{\frac{p-1}{2}}\frac{4^k}{k} \pmod{p}.  \label{5pr-11}
	\end{aligned}
\end{align}
Recall that the Granville congruence \cite{wt-12}:
\begin{align*}
		\sum_{k=1}^{p-1}\frac{x^k}{k}\equiv 	\frac{1-x^p-(1-x)^p}{p} \pmod{p}.
\end{align*}
It follows that
\begin{align}
	\begin{aligned}
		\sum_{k=1}^{\frac{p-1}{2}}\frac{4^k}{2k-1}+\sum_{k=1}^{\frac{p-1}{2}}\frac{4^k}{k}=2\sum_{k=1}^{p-1}\frac{2^k}{k} \equiv \frac{2-2^p}{p} \pmod{p}. \label{5pr-12}
	\end{aligned}
\end{align}
Combining   \eqref{5pr-8}--\eqref{5pr-12} gives
\begin{align}
	\begin{aligned}
		\sum_{k=1}^{\frac{p-3}{2}}F^*(p,k)+\sum_{k=\frac{p+1}{2}}^{p-2}F^*(p,k)\equiv -78+60q_p(2) \pmod{p}.  \label{5pr-13}
	\end{aligned}
\end{align}

Using \eqref{wl0-1} and \eqref{wl0-4}, we have
\begin{align}
	\begin{aligned}
		{2p-1\choose (p-1)/2}&\equiv (-1)^{\frac{p-1}{2}}\left(1-2pH_{\frac{p-1}{2}}\right)\equiv (-1)^{\frac{p-1}{2}}\left(1+4pq_p(2)\right)  \pmod{p^2}, \label{5pr-13-1}
	\end{aligned}
\end{align}
and
\begin{align}
	\begin{aligned}
		{(7p-1)/2\choose p}&\equiv 3+3pH_{\frac{p-1}{2}}+2pH_p\equiv 3-6pq_p(2)  \pmod{p^2}. \label{5pr-13-2}
	\end{aligned}
\end{align}
It is intutive
\begin{align}
	\begin{aligned}
	{4p\choose 2p}/{4p\choose p}\equiv {3p\choose 2p}/{2p\choose p}\equiv \frac{3}{2} \pmod{p^2}. \label{5pr-13-3}
	\end{aligned}
\end{align}
Utilizing  \eqref{4pr-19}, \eqref{5pr-6},  \eqref{5pr-13-1}--\eqref{5pr-13-3}, and Fermat's Little Theorem yields
\begin{align}
	\begin{aligned}
		F^*(p,p-1)=\frac{(3p+2)(p-1)}{4(p+1)^2p}{3p\choose p}{2p\choose p}\equiv \frac{9}{2}-\frac{p}{3} \pmod{p}, \label{5pr-14}
	\end{aligned}
\end{align}
and
\begin{align}
	\begin{aligned}
		F^*(p,\frac{p-1}{2})=\frac{(-1)^{\frac{p-1}{2}}(9p+1)(p-1)}{2^{p-2}3(3p+1)^2(3p-1)p}
		\frac{{4p\choose 2p}{\frac{7p-1}{2}\choose p}{2p-1\choose \frac{p-1}{2}}^2}{{4p\choose p}{p-1\choose \frac{p-1}{2}}} \equiv 15+\frac{3}{p}+9q_p(2) \pmod{p}. \label{5pr-15}
	\end{aligned}
\end{align}
Substituting \eqref{5pr-13}--\eqref{5pr-15} into \eqref{5pr-5} and simplifying, we obtain
\begin{align*}
	\begin{aligned}
		&\sum_{k=1}^{p-1}\frac{(-1)^k(3p+3k-1)}{2^{2k}(p+k-1)(p+k)(2p+2k-1)}\frac{{2p+2k\choose p+k}{3p+k-1\choose p}{p-1\choose k}}{{2k\choose k}}\equiv -24q_p(2)  \pmod{p}.
	\end{aligned}
\end{align*}
This, together with \eqref{4pr-19} yields
\begin{align}
	\begin{aligned}
		\sum_{k=0}^{p-1}&\frac{(-1)^k(3p+3k-1)}{2^{2k}(p+k-1)(p+k)(2p+2k-1)}\frac{{2p+2k\choose p+k}{3p+k-1\choose p}{p-1\choose k}}{{2k\choose k}}\\
		&\equiv \frac{2(3p-1)}{3p(p-1)(2p-1)}{2p\choose p}{3p\choose p}-24q_p(2) \\
		&\equiv -\frac{4}{p}-24q_p(2) \pmod{p}. \label{5pr-16}
	\end{aligned}
\end{align}
Then, combining  \eqref{4pr-12}, \eqref{4pr-19}, \eqref{5pr-2},  \eqref{5pr-3},  and \eqref{5pr-16}, we arrive at
\begin{align*}
	\begin{aligned}
		\sum_{k=2}^{p-1}\frac{11k+3}{2^{6k}}{2k\choose k}{3k\choose k}\equiv \frac{3p+18p^2q_p(2)}{2^{6(p-1)}}-\frac{45}{8} \pmod{p^3}. 
		\end{aligned}
	\end{align*}
Therefore, we obtain
\begin{align*}
	\begin{aligned}
		\sum_{k=0}^{p-1}\frac{11k+3}{2^{6k}}{2k\choose k}{3k\choose k}=\frac{45}{8}+ \sum_{k=2}^{p-1}\frac{11k+3}{2^{6k}}{2k\choose k}{3k\choose k}\equiv 3p \pmod{p^3}, 
	\end{aligned}
\end{align*}
as claimed result.	This completes the proof. \qed

\section{Proof of the Theorem 1.4}	
	
On the one hand, summing up both sides of the WZ equation \eqref{5pr-1} for $k$ form $0$ to $p-1$, and for $n$ from $2$ to $\frac{p-1}{2}$, we derive that
\begin{align*}
	\sum_{k=0}^{p-1}F^{'}(\frac{p+1}{2},k)-\sum_{k=0}^{p-1}F^{'}(2,k)=-\sum_{n=2}^{\frac{p-1}{2}}G^{'}(n,0), 
\end{align*}
with some arrangement and calculation, this is equivalent to	
\begin{align}
	\begin{aligned}
	\frac{{p-1\choose \frac{p-1}{2}}}{2^{3p+1}}&\sum_{k=0}^{\frac{p-1}{2}}\frac{(-1)^k(3p+1+6k)(p+1)p}{4^k(p-1+2k)(p+1+2k)^2}{p-1+2k\choose \frac{p-1+2k}{2}}{\frac{3p+1+2k}{2}\choose p+k}{\frac{p-1}{2}\choose k}/{2k\choose k}\\
	&=\frac{15}{128}-\sum_{k=2}^{\frac{p-1}{2}}\frac{22k^2-3k-3}{2^{6k+3}k(k-1)(2k-1)}{2k\choose k}^2{3k\choose k}. \label{6pr-1}
\end{aligned}	
\end{align}	
On the other hand, summing both sides of \eqref{5pr-2-1} over $n$ from $2$ to $\frac{p-1}{2}$, we obtain
\begin{align}
	\begin{aligned}
		\sum_{k=2}^{\frac{p-1}{2}}&\frac{11k+3}{2^{6k}}{2k\choose k}^2{3k\choose k}-48\sum_{k=2}^{\frac{p-1}{2}}\frac{22k^3-3k-3}{2^{6k+3}3k(k-1)(2k-1)}{2k\choose k}^2{3k\choose k} \\
		&=\frac{3(9p^2-1)}{2^{3p-2}(p-1)}{\frac{3p-3}{2}\choose \frac{p-1}{2}}{p-1\choose \frac{p-1}{2}}^2-\frac{45}{4}.  \label{6pr-2}
	\end{aligned}	
\end{align}		
	
For $1\leq k\leq \frac{p-1}{2}$, we have	
	\begin{align}
		\begin{aligned}
			{2p-k\choose \frac{p-1}{2}}&=\frac{(p+p-k)\cdots(p+1)(p-1)\cdots(p-\frac{p-1}{2})}{(p+\frac{p+1-k}{2})\cdots (p+1)(p-1)!}\\
			&\equiv (-1)^{\frac{p-1}{2}}\frac{2(p-k){\frac{p-1}{2}\choose k}}{(p+1-2k){p-1\choose k}} \pmod{p}.  \label{6pr-3}
		\end{aligned}	
	\end{align}			
A result due to Sun \cite[(4.4)]{wt-10} says: for $p=2n+1$ be an odd prime,	
\begin{align}
	(-4)^k{n\choose k}\equiv {2k\choose k} \left(1-p\sum_{j=1}^{k}\frac{1}{2j-1}\right)  \pmod{p^2}. 	  \label{6pr-4}
\end{align}	
Hence, 	using  \eqref{3pr-5}, \eqref{3pr-5-0}, \eqref{3pr-10}, \eqref{4pr-2-1}, \eqref{5pr-12} \eqref{6pr-3}, \eqref{6pr-4}, we obtain
\begin{align}
	\begin{aligned}
		\sum_{k=2}^{\frac{p-1}{2}}&\frac{(-4)^k(3p+2-3k)(2p+1-k)k}{(p-k)(p-2k)(p+1-k)^2}\frac{{\frac{p-1}{2}\choose k}{2p-2k\choose p-k}{2p-k\choose \frac{p-1}{2}}}{{p-1-2k\choose \frac{p-1-2k}{2}}}\\
		&\equiv -(-1)^{\frac{p-1}{2}}\frac{4p}{{p-1\choose \frac{p-1}{2}}}\left(\sum_{k=1}^{\frac{p-1}{2}}\frac{4^k}{k}+\sum_{k=1}^{\frac{p-1}{2}}\frac{4^k}{2k-1}+4\right)\\
		&\equiv -16p+16pq_p(2) \pmod{p^2}.   \label{6pr-5}
	\end{aligned}	
\end{align}			
Moreover, utilizing \eqref{3pr-10}, \eqref{3pr-12}, \eqref{3pr-13}, \eqref{5pr-13-1},  we arrive at 
\begin{align}
	\begin{aligned}
	\frac{-16(3p-1)}{(p-1)(2p-1)}{2p-1\choose p-1}{2p-1\choose \frac{p-1}{2}}/{p-1\choose \frac{p-1}{2}}\equiv 16+32pq_p(2) \pmod{p^2},  \label{6pr-6}
	\end{aligned}	
\end{align}	
\begin{align}
	\begin{aligned}
		\frac{(3p+1)(9p^2-1)}{2(p^2-1)(p+1)}{\frac{3p-3}{2}\choose \frac{p-1}{2}}{p-1\choose \frac{p-1}{2}}\equiv -(-1)^{\frac{p-1}{2}}\left(p+5p^2\right) \pmod{p^3},  \label{6pr-7}
	\end{aligned}	
\end{align}	
and
\begin{align}
	\begin{aligned}
		\frac{3(9p^2-1)}{2^{3p-2}(p-1)}{\frac{3p-3}{2}\choose \frac{p-1}{2}}{p-1\choose \frac{p-1}{2}}^2\equiv -3\left(p+4p^2-p^2q_p(2)\right) \pmod{p^3}.  \label{6pr-8}
	\end{aligned}	
\end{align}	
Performing the relation $k\to \frac{k+1-k}{2}$, and with the help of \eqref{6pr-5} and \eqref{6pr-6} we have
\begin{align*}
	\begin{aligned}
		\sum_{k=1}^{\frac{p-1}{2}}&\frac{(-1)^k(3p+1+6k)(p+1)p}{4^k(p-1+2k)(p+1+2k)^2}{p-1+2k\choose \frac{p-1+2k}{2}}{\frac{3p+1+2k}{2}\choose p+k}{\frac{p-1}{2}\choose k}/{2k\choose k}\\
		&=-(-1)^{\frac{p-1}{2}}\frac{p}{2^{p+3}}\sum_{k=1}^{\frac{p-1}{2}}\frac{(-4)^k(3p+2-3k)(2p+1-k)k}{(p-k)(p-2k)(p+1-k)^2}\frac{{\frac{p-1}{2}\choose k}{2p-2k\choose p-k}{2p-k\choose \frac{p-1}{2}}}{{p-1-2k\choose \frac{p-1-2k}{2}}}\\
		&\equiv -(-1)^{\frac{p-1}{2}}\left(p-p^2+2p^2q_p(2)\right) \pmod{p^3}.
	\end{aligned}	
\end{align*}		
This, together with \eqref{6pr-1}, \eqref{6pr-7} yields
 \begin{align}
 	\begin{aligned}
 	\sum_{k=2}^{\frac{p-1}{2}}\frac{22k^2-3k-3}{2^{6k+3}k(k-1)(2k-1)}{2k\choose k}^2{3k\choose k}
 	\equiv \frac{p}{8}+\frac{p^2}{4}+\frac{15}{128} \pmod{p^3}. \label{6pr-9}
 	\end{aligned}	
 \end{align}			
Substituting \eqref{6pr-8} and \eqref{6pr-9} into \eqref{6pr-2} simplifying, we obtain
 \begin{align*}
	\begin{aligned}
	\sum_{k=2}^{\frac{p-1}{2}}\frac{11k+3}{2^{6k}}{2k\choose k}^2{3k\choose k}\equiv -\frac{45}{8}+3p+3p^2q_p(2) \pmod{p^3}.
	\end{aligned}	
\end{align*}	
It follows that
 \begin{align*}
	\begin{aligned}
		\sum_{k=0}^{\frac{p-1}{2}}&\frac{11k+3}{2^{6k}}{2k\choose k}^2{3k\choose k}=\frac{45}{8}+\sum_{k=2}^{\frac{p-1}{2}}\frac{11k+3}{2^{6k}}{2k\choose k}^2{3k\choose k}
		\equiv 3p+3p^2q_p(2) \pmod{p^3}.
	\end{aligned}	
\end{align*}		
Now, the proof of Theorem $1.4$ is complete. \qed


\begin{thebibliography}{99}
		\small \setlength{\itemsep}{-.8mm}
		
		\bibitem{wt-0-1}Ramanujan S.: Modular equations and approximation to $\pi$, Quart. J. Math. \textbf{45} (1914),pp. 350--372.
		
		\bibitem{wt-0-2}Guillera J.: A method for proving Ramanujan's series for $1/\pi$, Ramanujan J. \textbf{52} (2020), pp. 421--431.
		
		\bibitem{wt-0-3}Van Hamme L.: Some conjectures concerning partial sums of generalized hypergeometric series, in p-Adic Functional Analysis (Nijmegen, 1996), Lecture Notes in Pure and Appl. Math. \textbf{192}, Dekker, New York, (1997), pp. 223--236.
		
		
		\bibitem{wt-0-4}Wilf H. S. and Zeilberger D.: Rational functions certify combinatorial identities, J. Am. Math.Soc. \textbf{3}(1) (1990), pp. 147--158.
		
		\bibitem{wt-0-5}Gessel I. M.: Finding identities with the WZ method, J. Symb. Comput. \textbf{20}(5) (1995), pp. 537--566.
		
		\bibitem{wt-3}Feng L.-Q  and Hou Q.-H.: Finding congruences with WZ method, Ramanujan J. \textbf{68}(3) (2025), pp. 64.  
		
		
		\bibitem{wt-4}Feng L.-Q  and Hou Q.-H.: Some supercongruences of hyperheometric sums and the WZ method, J. Differ. Equ. Appl. \textbf{31} (2025), pp. 1--26.
		
		\bibitem{wt-0-6}Guo V. J. W. and Liu J.-C.: Some congruences related to a congruence of Van Hamme, Integra Transforms and Special Functions, \textbf{31}(3) (2020), pp. 221-231. 
		
		
		\bibitem{wt-0-7}Mao G.-S.: Proof of some supercongruences via the Wilf-Zeilberger method, J. Differ. Equ. Appl. \textbf{26} (2020), pp. 1494--1513.
		
		\bibitem{wt-0-8}Sun Z.-W.: A refinement of a congruence result by van Hamme and Mortenson, Ill. J. Math. \textbf{56} (2012), pp. 967--979.
		
		\bibitem{wt-0-9}Zudilin, W.: Ramanujan-type supercongruences, J. Number Theory \textbf{129} (2009), pp. 1848–1857.
		
		\bibitem{wt-0-10}Wang C., Hu D.-W.: Proof some supercongruences concerning truncated hypergeometric series, Rev. Real Acad. Clenc. Exactas Fis. Nat.Ser. A-Mat., \textbf{117}(3) (2023), 99.
		
		
		
		
		\bibitem{wt-1}Guillera J. and  Zudilin W., `Divergent' Ramanujan-type supercongruences, Proc. Am. Math. Soc. \textbf{140}(3) (2012), pp. 765--777.
		
		\bibitem{wt-2}Bailey W. N., Generalized Hypergeometric Series, Cambridge University Press, Cambridge, (1935).
	
	
		\bibitem{wt-5}Wolstenholme J.: On certain properties of prime numbers. Q. J. Appl. Math. \textbf{5} (1862), pp.	35--39.
		
	
		\bibitem{wt-6}Mao G.-S.: Congruences involving Franel numbers and Apéry-like numbers,  temporarily on Researchgate (2023), preprint.	
		
	\bibitem{wt-7}Sun Z.-H.:  Congruences concerning Bernoulli numbers and Bernoulli polynomials, Discrete. Appl. Math. \textbf{105}  (2000), pp. 193--223	
	
		\bibitem{wt-8}Mao G.-S.: On some super-congruences for the coefficients of analytic solutions of certain differential equations, Indian J Pure Appl Math. \textbf{57} (2026), pp. 409--418.
		
		\bibitem{wt-9}Morely E.: Note on the congruence $2^{4n}\equiv (-1)^n(2n)!/(n!)^2$, where $2n+1$ is a prime, Ann. Math. \textbf{9} (1895), pp. 168--170.
		
	\bibitem{wt-10}Sun Z.-W: A new series for $\pi^3$ and related congruences. Internat. J. Math. \textbf{26}(8) (2015), 1550055 (23 pages). 
	
		
		\bibitem{wt-11}Chu W.-C.: Inversion techniques and combinatorial identities: A uniffed treatment for the ${}_7F_6$-series identities, Collect. Math. \textbf{45}, (1994), pp. 13--43.
		
		\bibitem{wt-12}Granville A.: The square of the Fermat quotient. Integers: Electronic Journal of Combi-natorial Number Theory \textbf{4}(A22) (2004), pp. 1--3.
		
		
		
		
		
	\end{thebibliography}
\end{document}